\documentclass[12pt,oneside,notitlepage]{amsart}
\usepackage{graphicx}
\newtheorem{thm}{Theorem}[section]
\newtheorem{cor}[thm]{Corollary}
\newtheorem{lem}[thm]{Lemma}
\newtheorem{prop}[thm]{Proposition}
\theoremstyle{definition}
\newtheorem{defn}[thm]{Definition}
\newtheorem{ex}[thm]{Example}
\theoremstyle{remark}
\newtheorem{rem}[thm]{Remark}
\numberwithin{equation}{section}
\newcommand{\norm}[1]{\left\Vert#1\right\Vert}

\newcommand{\set}[1]{\left\{#1\right\}}

\newcommand{\To}{\longrightarrow}

\newcommand{\A}{\mathcal{A}}
\newcommand{\Aut}{\textrm{Aut}}
\newcommand{\Herm}{\textrm{Herm}}
\newcommand{\tr}{\textrm{tr}}

\newcommand{\dist}{\textrm{dist}}
\newcommand{\CP}{\mathbb{CP}}

\newcommand{\R}{\mathbb{R}}
\newcommand{\C}{\mathbb{C}}

\newcommand{\M}{\mathcal{M}}
\newcommand{\K}{\mathcal{K}}
\newcommand{\B}{\mathcal{B}}
\newcommand{\T}{\mathcal{T}}
\newcommand{\pp}{\mathfrak{p}}

\newcommand{\Hilb}{\mathrm{Hilb}}
\newcommand{\FS}{\mathrm{FS}}

\begin{document}

\title[Balanced metrics on non-simple bundles]{Numerical algorithm for finding balanced metrics on vector bundles II: the non-simple case}%

\author[R. Seyyedali]{Reza Seyyedali}
\address {School of Mathematics, Institute for Research in Fundamental Sciences (IPM), P.O. Box: 19395-5746, Tehran, Iran}
\email{rseyyedali@ipm.ir}
\address{Department of Mathematics, New Uzbekistan University, Tashkent,  Uzbekistan}
\email{r.seyyedali@newuu.uz}


\begin{abstract}

In \cite{D3}, Donaldson defines a dynamical system on the space of
Fubini-Study metrics on a polarized compact K\"ahler manifold.
Sano proved that if there exists a balanced metric for the
polarization, then this dynamical system always converges to the
balanced metric provided that the automorphism group is discrete (\cite{S}).
 In \cite{DKLR}, Douglas, et.\ al., conjecture that the same holds in the case of vector bundles. In \cite{Se}, we proved their conjecture under the assumption that the bundle is simple. In this note, we remove the simplicity assumption. The main tool is Thompson's metric on the space of Hermitian inner products on the space of holomorphic sections, with respect to which the generalized $T$-operator is non-expanding. 

\end{abstract}
\maketitle
\pagestyle{plain} \pagenumbering{arabic}

\section{Introduction}

The notion of stability is a  central notion in geometric invariant theory (GIT) which was invented to construct quotient
spaces of projective varieties.  An important example of stability notion is the  Mumford-Takemoto slope stability for holomorphic vector bundles.  There is also another notion of stability introduced by  Gieseker 
which is more in the realm of geometric invariant theory. A celebrated theorem of  Donaldson, Uhlenbeck and Yau
(\cite{D1}, \cite{D2}, \cite{UY}) establishes that over a compact K\"ahler manifold,
a holomorphic vector bundle is poly-stable if and only if it admits a
Hermitian-Einstein metric. This is known as the Hitchin-Kobayashi
correspondence.  In \cite{W}, Wang proves a conjecture of Donaldson that states a vector bundle $E$ over a polarized manifold $(X,L)$ is Gieseker stable if and only if $E\otimes L^{k}$ admits a balanced metric for $k\gg 0$. He further proves that on a stable vector bundle, the sequence of rescaled balanced metrics converges to the Hermitian-Einstein metric (\cite{W2}). 
In \cite{D3}, Donaldson defines a dynamical system on the space of
Fubini-Study metrics on a polarized compact K\"ahler manifold. Sano proved that this dynamical system converges whenever a balanced metric exists and the automorphism group is discrete (\cite{S}). In \cite[Proposition 4]{D3}, Donaldson also proves the convergence of the analogous map $T_{\nu}$, which integrates against a fixed measure $\nu$, under assumptions on $\nu$ which imply that the $\nu$-balanced metric is unique up to scaling (see Remark \ref{genrem}). In \cite[Conjecture 3.5]{DKLR}, Douglas, et.\ al., conjecture that the same holds in the case of vector bundles. They assume that the bundle is irreducible, and they remark that the conjecture may require additional assumptions, such as the discreteness of the automorphism group. This conjecture is proved in \cite{Se} under the assumption that the bundle is simple. The goal of this note is to prove their conjecture in full generality. In particular, no assumption on the automorphism group of $E$ is needed.

Let $(X,\omega)$ be a compact K\"ahler manifold of dimension $m$ and $E$
be a very ample holomorphic vector bundle of rank $r$ on $X$. Let $h$ be a
Hermitian metric on $E$. We can define an $L^2$-inner product on
$H^{0}(X,E)$ by

$$\langle s,t\rangle= \frac{N}{Vr}\int_{X} h(s,t)\frac{\omega^m}{m!}.$$
Let $s_{1},...,s_{N}$ be an orthonormal basis for $H^{0}(X,E)$
with respect to this $L^2$-inner product. The Bergman kernel of
$h$ is defined by $$B(h)=\sum s_{i}\otimes s_{i}^{*_{h}}.$$ Note
that $B(h)$ does not depend on the choice of the orthonormal basis
$s_{1},...,s_{N}$.

A metric $h$ is called balanced if $B(h)$ is a constant multiple
of the identity. By the work of Wang and Phong--Sturm (\cite{W}, \cite{PS}),
the existence of balanced metrics is closely related
to the stability of the vector bundle $E$. Indeed, by \cite[Theorem 2]{PS} and the Kempf--Ness theorem (\cite{KN}), $E$ admits balanced metrics if and only if
the Gieseker point of $E$ is poly-stable (c.f. \cite[Theorem 1.1]{W}, whose proof treats the stable case). Note that the uniqueness statement of \cite[Theorem 1.1]{W} does not hold for strictly poly-stable bundles, such as those in Example \ref{excp} with $k\geq2$ (see Proposition \ref{nsprop1}). Note also that this is a GIT condition on the bundle $E$ itself, while Gieseker stability is a condition on $E\otimes L^{k}$ for all $k\gg 0$. Let $\K$ and $\M$ be the spaces of
Hermitian metrics on $E$ and Hermitian inner products on
$H^{0}(X,E)$ respectively. Following Donaldson (\cite{D3}), one
can define the following maps

\begin{itemize}

\item Define $$\textrm{Hilb}: \K \rightarrow \M$$ by

 $$\langle s,t \rangle _{\Hilb(h)}= \frac{N}{Vr} \int
\langle s(x),t(x) \rangle_{h}\frac{\omega^m}{m!},$$where
$N=\dim(H^{0}(X,E)) $ and $V=\textrm{  Vol}(X,\omega)$. Note that
$\textrm{Hilb}$ only depends on the volume form
$\frac{\omega^{m}}{m!}$.

\item For the metric $H$ in $\M$, the metric $\FS(H) \in \K$ is the unique metric
on $E$ such that  $\sum s_{i}\otimes s_{i}^{*_{\FS(H)}}=I$, where
$s_{1},...,s_{N}$ is an orthonormal basis for $H^{0}(X,E)$ with
respect to $H$. This gives the map $\FS:\M \rightarrow \K$.

\item Define a map $$T:\M \rightarrow \M$$ by $T(H)=\Hilb
\circ \FS(H).$ Notice that this map $T$ is called generalized
$T$-operator in \cite{DKLR}.

\end{itemize}

If $H\in \M$ is a fixed point of the map $T$, then $\FS(H)$ is balanced and $H=\Hilb(\FS(H))$. Conversely, if $h$ is balanced, then $\FS(\Hilb(h))=h$ and $\Hilb(h)$ is a fixed point of $T$. Therefore we also call the fixed points of $T$ balanced metrics.

The main theorem of this paper is the following
\begin{thm}\label{mainthm}

Suppose that $E$ admits a balanced metric $H_{\infty}$. Then for
any $H_{0}\in \M$, the sequence $T^n(H_{0})$ converges to
$\varphi^*H_{\infty}$ for some $\varphi \in \Aut(E)$ whose action on $H^{0}(X,E)$ is self-adjoint and positive with respect to $H_{\infty}$.

\end{thm}

In \cite{Se}, we proved Theorem \ref{mainthm} for simple vector bundles. 
The proof there follows  Sano's argument in \cite{S} with the necessary modifications for the bundle case. The main tool in the proof is the functional $\widetilde{Z}$ which
is introduced by Wang (\cite{W}) and Phong, Sturm (\cite{PS})
to study the existence and uniqueness of balanced metrics on
holomorphic vector bundles. The key property of this functional is
that it is convex along geodesics and its critical points are balanced metrics. By the Kempf--Ness theorem, one can show that this functional is proper modulo scaling if $E$ is simple and admits a balanced metric. Therefore, if $E$ is simple, the sequence $H_n=T^n(H_0)$ is bounded modulo scaling, since $\widetilde{Z}(H_n)$ is decreasing and bounded from below. 
However, if $E$ is not simple, then the space of balanced metrics modulo scaling is non-compact, and therefore $\widetilde{Z}$ is neither strictly convex nor proper. Hence, the boundedness of $\widetilde{Z}(H_n)$ does not imply the boundedness of $H_n.$ In order to overcome this issue, we use a distance function $\dist_T$ on $\M$ introduced by Thompson in \cite{T} (see also \cite{N} and \cite{LN}). The advantage of this distance function is that the map $T$ is non-expanding with respect to it and therefore, the sequence $H_n$ is bounded (Proposition \ref{nsprop2}). Once we prove this proposition, it suffices to show that  subsequential limits of $H_n$ are balanced and that the whole sequence converges to a unique limit. One distinguishing feature of the current argument is that we do not use the functional $\widetilde{Z}$ in order to prove that the sequence $\{H_n\}$ converges to a balanced metric. Indeed the only ingredients used in the proof of convergence are the existence of balanced metrics, order preservation and homogeneity of $T$, and the trace identity $\tr(T(H)H^{-1})=N$. The only place that we use the functional $\widetilde{Z}$ is the proof of Proposition \ref{nsprop1}, which identifies the limit.
Note that the argument is close in spirit to some results outside complex geometry. By Lemma \ref{partialorder}, $T(H)(s,t)=H(A_{H}s,t)$, where
$$A_{H}=\frac{N}{Vr}\int_{X}P_{x}\,\frac{\omega^{m}}{m!}$$
and $P_{x}$ is the $H$-orthogonal projection onto the $H$-orthogonal complement of $K_x=\set{s\in H^{0}(X,E)\mid s(x)=0}$. Hence $H$ is balanced if and only if $A_{H}=I$. If the volume form is replaced by a finite sum of point masses, this says that finitely many $r$-dimensional subspaces are in isotropic position. For $r=1$, the map $T$ is then Tyler's fixed point iteration for the $M$-estimator of scatter (\cite{Ty}). In particular, Remark \ref{genrem}(2) shows that Tyler's iteration converges whenever the data points are not contained in a hyperplane and the $M$-estimator exists, also when it is not unique up to scaling. For $r>1$, the condition $A_{H}=I$ can be regarded as a complex analogue of the condition defining geometric Brascamp-Lieb data (\cite{BCCT}), and the problem of moving given data to such a position is studied by operator scaling (\cite{Gu}, \cite{GGOW}). The outline of the paper is as follows.

In section $2,$ we  some basic definitions and properties of $T$ and $\widetilde{Z}.$ We also give a description of the set of balanced metrics in $\M$ (Proposition \ref{nsprop1}).
In section $3,$ we prove Theorem \ref{mainthm}. We first introduce the distance function $\dist$ on $\M$ (Definition \ref{distance}) and show that it coincides 
with Thompson's metric $\dist_T$ \eqref{thompsondistance}. This together with a simple observation that $T$ preserves the partial order \eqref{partialorderdef}
on $\M$ (Lemma \ref{partialorder}) implies that $T$ is non-expanding (Proposition \ref{nsprop2}). At the end of section $3$, we discuss some generalizations of Theorem \ref{mainthm} (Remark \ref{genrem}).
In section $4,$ we sketch an argument to show that the convergence is exponential, and we compute the rate for $\mathcal{O}(d)\otimes\C^{k}$ over $\CP^{m}$ (Example \ref{excp}). We also show that the asymptotic formula of \cite[Theorem 1.2(1)]{IP} for the rate holds for polystable bundles (Proposition \ref{proppoly}), and we give an example in which the rate is determined by the homomorphisms between the stable summands (Example \ref{exP1P1}).
In the appendix, we compute the variation of $\widetilde{Z}$ along the orbits of $\Aut(E)$ (Proposition \ref{propinv}), and we obtain an obstruction to the existence of balanced metrics (Corollary \ref{cordirect}).

\medskip
\noindent\textbf{Acknowledgements:} I used the AI assistant Claude (Opus 5.5, Anthropic) in the preparation of this paper. The references to work outside complex geometry were suggested by Claude. Claude was also used to expedite routine computations in the examples and proofs of Section \ref{secrate} and Appendix \ref{appinv}, to check numerically some cases of the rates of convergence computed in the examples of Section \ref{secrate}, and to review a draft of the paper. It also assisted in identifying and correcting minor errors and in filling minor gaps. I emphasize that the main results and their proofs, the overall structure of the paper and the vast majority of the examples are my own work, with some help from Claude in expediting calculations. I have checked all of the mathematics and the references, and I take full responsibility for the content of the paper.

\section{Balanced metrics on vector bundles}

As before, let $(X,\omega)$ be a compact K\"ahler manifold of dimension $m$ and $E$ be a
very ample holomorphic vector bundle on $X$. Using global sections
of $E$, we can map $X$ into $G(r, H^{0}(X,E)^*)$. Indeed, for any
$x \in X$, we have the evaluation map $$s \in H^{0}(X,E)\rightarrow
s(x) \in E_{x}. $$ Since $E$ is globally
generated, this map is a surjection. So its dual is an inclusion
of $E_{x}^* \hookrightarrow H^{0}(X,E)^*$, which determines a
$r$-dimensional subspace of $H^{0}(X,E)^*$. Therefore we get an
embedding $$\iota: X\hookrightarrow G(r, H^{0}(X,E)^*).$$ Clearly we
have $\iota^*U_{r}=E^*$, where $U_{r}$ is the tautological vector
bundle on $G(r, H^{0}(X,E)^*)$, i.e. at any $r$-plane in $G(r,
H^{0}(X,E)^*)$, the fibre of $U_{r}$ is exactly that $r$-plane. A
choice of basis for $H^{0}(X,E)$ gives an isomorphism between
$G(r, H^{0}(X,E)^*)$ and the standard $G(r,N)$, where $N= \dim
H^{0}(X,E)$. We have the standard Fubini-Study hermitian metric on
$U_{r}$, so we can pull it back to $E$ and get a hermitian metric
on $E$.  Using $\iota^*h_{\FS}$ and $\omega$, we get an $L^{2}$ inner
product on $H^{0}(X,E)$. The embedding is called balanced if
$\int_{X} \, \langle s_{i}, s_{j}\rangle \,
\frac{\omega^{m}}{m!}=C \delta_{ij} $. We can formulate this
definition in terms of maps $\textrm{Hilb}$ and $\FS$.

 \begin{defn}

An $\omega$-balanced metric on $E$ is a pair $(h_b,H_b)$ such that
$$\Hilb(h_b)=H_b,\qquad \FS(H_b)=h_b.$$

 \end{defn}

 Fixing a nonzero element $\Theta \in \bigwedge^N H^{0}(X,E)$,
 we can define the determinant of any element in $\M$. Thus we
 can define a map $$\log \det: \M\rightarrow \mathbb{R}$$
 A different choice of $\Theta$ only changes this map by an additive
 constant.
 Also, we define a functional $I:\K\rightarrow \mathbb{R}$
 unique up to an additive constant. Fix a background metric
 $h_{0}$ and consider a path $h_{t}=e^{\phi_{t}}h_{0}$ in $\K,$
 \begin{equation}\label{eq1}\frac{dI}{dt}=\int_{X} \tr(\dot{\phi}) \textrm{   dVol}_{\omega}\end{equation}
 This functional is part of  Donaldson's functional.
 We define: \begin{equation}\label{eq2} Z=-I \circ \FS:\M\rightarrow \mathbb{R}\end{equation}

 We have the following scaling identities:
 $$\textrm{Hilb}(e^{\alpha}h)=e^{\alpha} \textrm{Hilb}(h),$$
 $$ \FS(e^{\alpha}H)=e^{\alpha} \FS(H),$$
 $$I(e^{\alpha}h)=I(h)+\alpha rV  ,$$
 where $\alpha$ is a real number.
Following Donaldson, we define:
\begin{equation}\label{eq3}\widetilde{Z}=Z+\frac{rV}{N}\log\det.\end{equation} So
$\widetilde{Z}$ is invariant under constant scaling of the metric.

This functional $Z$ is studied by Wang in \cite{W} and Phong and
Sturm in \cite{PS}. They consider this as a functional on
$SL(N)/SU(N)$. In order to see this, we observe that there is a
correspondence between $\M$ and $GL(N)/U(N)$. Fix an element
$H_{0} \in \M$ and an orthonormal basis $s_{1},...,s_{N}$ for
$H^{0}(X,E)$ with respect to $H_{0}$. Now for any $H \in \M$, the matrix $[H(s_{i},s_{j})]$ is positive definite Hermitian. Hence it can be written as $gg^{*}$ for some $g\in GL(N)$, and $g$ is unique up to right multiplication by an element of $U(N)$. Notice that change of the orthonormal basis only changes $g$ by multiplication by elements of $U(N)$. So we get a well-defined element $gU(N)$ of $GL(N)/U(N)$. The subset $$\M_{0}=\{H \in \M \mid \det [H(s_{i},
s_{j})]=1\}$$ corresponds to $SL(N)/SU(N)$.

We recall the definition of the Gieseker point of the bundle $E$ (c.f. \cite{G}).
We have a natural map
$$\T(E): \bigwedge^{ r}H^{0}(X,E)\rightarrow H^{0}(X,\det(E))$$ which
for any $s_{1},..., s_{r}$ in $H^{0}(X,E)$ is defined by
$$\T(E)(s_{1}\wedge...\wedge s_{r})(x)=s_{1}(x)\wedge...\wedge s_{r}(x).$$

Since $E$ is globally generated, $\T(E)$ is non-zero. We can view
$\T(E)$ as an element of
$\textrm{Hom}(\bigwedge^{r}H^{0}(X,E),H^{0}(X,\det(E)))$. This is
called the Gieseker point of $E$. Notice that fixing a basis for
$H^0(X,E)$ gives an isomorphism between $\bigwedge^{r}H^{0}(X,E)$
and $\bigwedge^r\mathbb{C}^N$. Hence, there is a natural action of
$GL(N)$ on
$\textrm{Hom}(\bigwedge^{r}H^{0}(X,E),H^{0}(X,\det(E)))$.
Phong-Sturm (\cite{PS}) and Wang (\cite{W}) prove that $Z$ is
convex along geodesics of $SL(N)/SU(N)$ and its critical points
are corresponding to balanced metrics on $E$. Phong and Sturm
prove the following

\begin{thm}(\cite [Theorem 2]{PS})\label{secondthm}
There exist a constant $c>0$, depending only on $E$ and $\omega$, and an $SU(N)$-invariant norm $\norm{\cdot}$ on the space
$$\textrm{Hom}\Big(\bigwedge^{r}H^{0}(X,E),H^{0}(X,\det(E))\Big)$$
such that for any $\sigma \in SL(N)$
$$Z(\sigma)-Z(I)= c\log \frac{\norm{\sigma\cdot\T(E)}^2}{\norm{\T(E)}^2},$$
where we identify $\sigma$ with the corresponding element of $\M_0$.
\end{thm}
In \cite{PS}, this identity is stated for the dual bundle, i.e. for $E\subseteq\C^N\times X$ with $E^*$ globally generated, with an explicit constant and norm.

\begin{rem}

In \cite{W}, Wang proves a slightly weaker version of Theorem
\ref{secondthm}. He proves that for any norm $||.||$ on the space
$$\textrm{Hom}\big(\bigwedge^{r}H^{0}(X,E),H^{0}(X,\det(E))\big),$$ there
exist positive constants $c$ and $c'$ such that

$$Z(\sigma) \geq c\log ||\sigma .
\T(E)||^2 +c'. $$

\end{rem}

\begin{thm}(\cite[Lemma 3.5]{W},  \cite[Lemma 2.2]{PS})\label{thirdthm}
The functional $Z$ is convex along geodesics of $\M$.

\end{thm}

The Kempf--Ness theorem (\cite{KN}) shows that $Z|_{\M_{0}}$ is proper and
bounded from below if $\T(E)$ is stable under the action of
$SL(N)$.

The first part of the following theorem is a consequence of Theorem \ref{thirdthm} and the fact that balanced metrics are critical points of $\widetilde{Z}$ (Lemma \ref{lemlin}(1) below). The second part follows from the Kempf--Ness theorem, since the stabilizer of the Gieseker point in $SL(N)$ is finite if $E$ is simple; it is used only in the account of the simple case. Note that the functional  $\widetilde{Z}$ is invariant under the scaling of a
metric by a positive real number.

\begin{thm}\label{fourththm}

Assume that $E$ admits balanced metrics. Then
$\widetilde{Z}|_{\M}$ is bounded from below, i.e. $\widetilde{Z}(H)\geq \widetilde{Z}(H_{b})$ for any $H
\in \M$ and any balanced metric $H_b$.
Moreover, $\widetilde{Z}|_{\M_0}$ is proper if $E$ is simple. 

\end{thm}
We recall the following lemma from \cite{Se}.

\begin{lem}\label{firstlemma}(\cite{Se}) For any $H \in \M$,
\begin{enumerate}

\item $\tr(T(H)H^{-1})=N.$
\item $Z(H)\geq Z(T(H)).$
 \item $\log\det (H) \geq\log\det(T(H)).$ Moreover, the equality holds if and only if $T(H)=H.$
\item $\widetilde{Z}(H)\geq \widetilde{Z}(T(H)).$
\end{enumerate}

\end{lem}

\begin{rem}
The proof that $T^{n}(H_{0})$ converges to a balanced metric uses only the existence of a balanced metric, the order preservation and homogeneity of $T$ (Lemma \ref{partialorder}), and parts $(1)$ and $(3)$ of Lemma \ref{firstlemma}. In particular, it does not use the functional $\widetilde{Z}$. The functional $\widetilde{Z}$ enters only through Proposition \ref{nsprop1}, which identifies the limit as $\varphi^{*}H_{\infty}$.
\end{rem}
A bundle $E$ is called simple if $\Aut(E)\simeq \mathbb{C}^*$.
If $E$ is a simple bundle  admitting balanced metrics, then the
balanced metric is unique up to a positive constant. 
In general the group $\Aut(E)$ acts on $\M$ and the set of balanced metrics $\B$ is the orbit of a balanced metric under this action. We prove that balanced metrics are unique up to the action of $\Aut(E)$ (Proposition \ref{nsprop1}).   

For $\varphi\in \textrm{Aut}(E)$, $h\in \K$ and $H\in \M$, define
$$\varphi^{*}h(v,w)=h(\varphi v,\varphi w), \,\,\,\,\,\,\,\,
\varphi^{*}H(s,t)=H(\varphi s,\varphi t),$$ where $\textrm{Aut}(E)$
acts on $H^{0}(X,E)$ by $s\mapsto \varphi\circ s$.

\begin{lem}\label{nslem4}
For any $\varphi\in \textrm{Aut}(E)$, $h\in \K$ and $H\in \M$, we have
$$\Hilb(\varphi^{*}h)=\varphi^{*}\Hilb(h), \quad
\FS(\varphi^{*}H)=\varphi^{*}\FS(H), \quad
T(\varphi^{*}H)=\varphi^{*}T(H).$$ In particular
$T^{n}(\varphi^{*}H_{0})=\varphi^{*}T^{n}(H_{0})$, and
$\varphi^{*}H$ is balanced whenever $H$ is balanced.
\end{lem}

\begin{proof}
The first identity is clear, since
$$\textrm{Hilb}(\varphi^{*}h)(s,t)=\frac{N}{Vr}\int_{X}h(\varphi s,\varphi
t)\frac{\omega^{m}}{m!}=\textrm{Hilb}(h)(\varphi s,\varphi t).$$ For the
second one, put $h=\FS(H)$ and let $s_{1},...,s_{N}$ be an orthonormal
basis with respect to $\varphi^{*}H$. Then $\varphi s_{1},...,\varphi
s_{N}$ is $H$-orthonormal. Hence the definition of $\FS(H)$ implies that for any
$v\in E_{x}$,
$$\varphi v=\sum_{i=1}^{N}h(\varphi v,\varphi s_{i}(x))\,\varphi
s_{i}(x).$$ Applying $\varphi^{-1}$, we get
$v=\sum\varphi^{*}h(v,s_{i}(x))\,s_{i}(x)$, i.e. $\sum s_{i}\otimes
s_{i}^{*_{\varphi^{*}h}}=I$. Therefore $\FS(\varphi^{*}H)=\varphi^{*}h$.
The third identity follows from the first two.
\end{proof}
The proof of the following proposition follows the argument of \cite{Se} in the simple case. For completeness, we give the proof here.
\begin{prop}\label{nsprop1}
Suppose that $E$ admits a balanced metric $H_{\infty}$. Then the set of
balanced metrics on $E$ is the orbit of $H_{\infty}$ under
$\Aut(E)$, i.e.
$$\B:=\set{H\in \M\mid T(H)=H}=\set{\varphi^{*}H_{\infty}\mid \varphi\in
\textrm{Aut}(E)}.$$ Moreover, for any balanced metric $H$ there exists
$\varphi\in \Aut(E)$, self-adjoint and positive with respect to
$H_{\infty}$, such that $H=\varphi^{*}H_{\infty}$.
\end{prop}

\begin{proof}
By Lemma \ref{nslem4}, $\varphi^{*}H_{\infty}$ is balanced for any
$\varphi\in \textrm{Aut}(E)$. Conversely, let $H$ be a balanced metric.
There exists a unique endomorphism $A$ of $H^{0}(X,E)$, self-adjoint
and positive with respect to $H_{\infty}$, such that
$H(s,s')=H_{\infty}(As,As')$. Write $A=e^{a}$, where $a$ is self-adjoint
with respect to $H_{\infty}$, and define
$$H_{t}(s,s')=H_{\infty}(e^{ta}s,e^{ta}s'),\,\,\,\,\,\,\,\, 0\leq
t\leq 1.$$ Then $H_{0}=H_{\infty}$, $H_{1}=H$ and $t\mapsto H_{t}$ is
a geodesic in $\M$. Put $f(t)=\widetilde{Z}(H_{t})$. Since $\log\det
H_{t}=\log\det H_{\infty}+2t\,\tr(a)$ is affine in $t$, Theorem
\ref{thirdthm} implies that $f$ is convex. Balanced metrics are critical
points of the convex functional $\widetilde{Z}$ (Lemma \ref{lemlin}(1)) and therefore minima
(c.f. Theorem \ref{fourththm}). Hence $f(0)=f(1)=\inf\widetilde{Z}$,
and convexity implies that $f$ is constant on $[0,1]$. Therefore
$$\frac{d^{2}}{dt^{2}}Z(H_{t})=\frac{d^{2}}{dt^{2}}f(t)=0,\,\,\,\,\,\,\,\, 0<t<1.$$
Let $\xi_a$ be the holomorphic vector field induced by $a$ on  $G(r,H^{0}(X,E)^{*}).$ By \cite[Lemma 3.5]{W}, we have 
$$\frac{d^{2}}{dt^{2}}Z(H_{t})=\int_X |\xi_a(x)|_{t}^2\frac{\omega^m}{m!},$$
where $|\xi_a(x)|_{t}$ is the Fubini-Study norm of the vector field $\xi_a$ at $\iota(x)$ with respect to $H_t$. The integrand is continuous and non-negative. Therefore, the vector field
$\xi_{a}$ vanishes on
$\iota(X)$ and hence the one parameter group $e^{ta}$ fixes every point of
$\iota(X)$. For $x\in X$, the point $\iota(x)=\iota_{x}(E_{x}^{*})$ is the
annihilator of $K_{x}=\set{s\in H^{0}(X,E)\mid s(x)=0}$. Therefore
$e^{ta}$ preserves $K_{x}$ for every $x\in X$ and $t\in\mathbb{R}$.
Note that $E$ is the quotient of the trivial bundle $H^{0}(X,E)\otimes
\mathcal{O}_{X}$ by the holomorphic subbundle $\mathcal{K}$ with fibres
$K_{x}$. Hence $e^{ta}$ descends to a holomorphic automorphism
$\varphi_{t}$ of $E$, and $(e^{ta}s)(x)=\varphi_{t}(s(x))$ for any $s\in
H^{0}(X,E)$. In other words, $e^{ta}$ is the action of $\varphi_{t}\in
\textrm{Aut}(E)$ on $H^{0}(X,E)$. Therefore
$$H(s,s')=H_{\infty}(e^{a}s,e^{a}s')=\varphi_{1}^{*}H_{\infty}(s,s'),$$
and $\varphi_{1}=e^{a}$ is self-adjoint and positive with respect to
$H_{\infty}$.
\end{proof}

\section{Proof of Theorem \ref{mainthm}}

 The main goal of this section is to prove Theorem \ref{mainthm}.  First, we recall the argument in the case of simple bundles. The first step in the simple case is to show that the sequence $H_n$ is bounded modulo scaling. The sequence $\widetilde{Z}(H_n)$ is bounded since the functional $\widetilde{Z}$ is bounded from below and the sequence  $\widetilde{Z}(H_n)$ is decreasing. On the other hand, the simplicity of the bundle implies that the functional $\widetilde{Z}|_{\M_0}$ is  proper. Hence the sequence $H_n$ is bounded modulo scaling. By Lemma \ref{firstlemma}, $\log\det(H_n)$ converges and therefore, the uniqueness of balanced metrics implies that the limit point of $\{H_n\}$ is unique. Therefore, $H_n$ converges to a balanced metric. In the general case, the proof fails for two main reasons. The most subtle difficulty arises from the fact that $\widetilde{Z}|_{\B}$ is constant, while $\B$ is non-compact even modulo scaling if $E$ is not simple. Therefore $\widetilde{Z}|_{\M_0}$ is neither proper nor strictly convex. This causes the main issue since the boundedness of $\widetilde{Z}(H_n)$ does not necessarily imply boundedness of $H_n$. Another issue is that the sequence $\{H_n\}$ could have more than one limit point and the sequence could potentially have different subsequential limits.
In order to overcome these issues, we show that the map $T$ is non-expanding with respect to the natural distance $\dist(H,H^{\prime})$ on $\M$ defined as follows. 

\begin{defn}\label{distance}
For $H,H'\in \M$, let $s_{1},\dots,s_{N}$ be an $H$-orthonormal basis and let $\mu_{1},
\dots,\mu_{N}$ be the eigenvalues of the positive definite matrix
$[H'(s_{i},s_{j})]$. Define
$$\dist(H,H')=\max_{1\leq i\leq N}|\log\mu_{i}|.$$
\end{defn}
Note that the eigenvalues $\mu_{i}$ do not depend on the choice of the $H$-orthonormal basis, since a change of such a basis conjugates the matrix $[H'(s_{i},s_{j})]$ by a unitary matrix. We show in Lemma \ref{distlem} that $\dist$ is a distance on $\M$.

\begin{prop}\label{nsprop2}
For any $H,H'\in \M$, we have
$$ \dist (T(H),T(H'))\leq \dist (H,H').$$
\end{prop}

With Proposition \ref{nsprop2} at hand, we conclude that the sequence $\{H_n\}$ is bounded. Then using the fact that $\log \det (H_n)$ is decreasing, we can show that all limit points of $H_n$ are balanced metrics. The remaining issue is to show the uniqueness of the limit points.
The main idea in the proof of Proposition \ref{nsprop2} comes from the work of Thompson \cite{T}. In \cite{T}, Thompson introduces a  metric on the interior of a cone
in a partially ordered vector space and proves that order preserving
homogeneous maps do not increase this distance. This is a basic tool of non-linear Perron-Frobenius theory (\cite{N}, \cite{LN}). It is also used in \cite{SH} for fixed point iterations on positive definite matices. We observe that the natural metric $\dist$ on $\M$ coincides with Thompson's metric.

We define a partial order on $\M$ by \begin{equation}\label{partialorderdef}H \leq H^{\prime} \,\,\,\,\textrm{if}  \,\,\,\,\, H(s,s) \leq  H^{\prime}(s,s),   \,\, \textrm{for all}\,\,\,\, s \in H^0(X,E). \end{equation} Note that it suffices to check the inequality only on the unit sphere in $H^0(X,E)$ with respect to $H$, which is compact. Hence, for any $H,H'\in \M$, there exists a positive number $\lambda$ such that $e^{-\lambda}H^{'} \leq H \leq e^{\lambda}H^{'}.$ Similarly, for $h,h'\in \K$ we write $h\leq h'$ if $h(v,v)\leq h'(v,v)$ for all $x\in X$ and $v\in E_{x}$. The maps $\FS$, $\Hilb$ and $T$ preserve these partial orders. 
\begin{lem}\label{partialorder}
Let $H,H'\in \M$. If $H\leq H'$, then $\FS(H)\leq \FS(H')$ and $T(H)\leq T(H')$. Moreover, $T(e^{\alpha}H)=e^{\alpha}T(H)$ for any $\alpha\in \R$, and
$$T(H)(s,t)=H(A_{H}s,t),\qquad A_{H}=\frac{N}{Vr}\int_{X}P_{x}\,\frac{\omega^{m}}{m!},$$
where $K_{x}=\set{s\in H^{0}(X,E)\mid s(x)=0}$ and $P_{x}$ is the $H$-orthogonal projection onto the $H$-orthogonal complement of $K_{x}$. In particular, $\tr(T(H)H^{-1})=\tr A_{H}=N$.
\end{lem}

\begin{proof}
By definition, $h=\FS(H)$ is the quotient metric induced by $H$ via the evaluation map $H^{0}(X,E)\rightarrow E_{x}$. Therefore, for $x\in X$ and $v\in E_{x}$, we have
$$|v|^{2}_{\FS(H)}=\min\set{H(s,s)\mid s\in H^{0}(X,E),\,\, s(x)=v}.$$ 
Indeed, $|v|^{2}_{\FS(H)}=H(\sigma_v,\sigma_v)$, where $\sigma_{v}=\sum_{i}h(v,s_{i}(x))\,s_{i}$ for $s_{1},\dots,s_{N}$ an $H$-orthonormal basis $s_{1},\dots,s_{N}$.  Hence $H\leq H'$ implies $\FS(H)\leq \FS(H')$. Since $\Hilb$ is clearly order preserving on the space of hermitian metrics on $E$, so is $T=\Hilb\circ \FS$. The homogeneity of $T$ follows from the scaling identities in Section 2. 

Finally, the map $v\mapsto\sigma_{v}$ is linear and $\sigma_{v}$ is $H$-orthogonal to $K_{x}$. Hence $\FS(H)(v,w)=H(\sigma_{v},\sigma_{w})$ by polarization. Moreover $\sigma_{s(x)}=P_{x}s$, since both are $H$-orthogonal to $K_{x}$ and $\sigma_{s(x)}-P_{x}s=(\sigma_{s(x)}-s)+(s-P_{x}s) \in K_x$.  Therefore
$$T(H)(s,t)=\frac{N}{Vr}\int_{X}H(P_{x}s,P_{x}t)\,\frac{\omega^{m}}{m!}=H(A_{H}s,t),$$
where we used $P_{x}^{2}=P_{x}$ and that $P_{x}$ is $H$-self-adjoint. Since $\tr P_{x}=\dim K_{x}^{\perp}=r$, we get $\tr A_{H}=N$.
\end{proof}

We now compute the first variation of $\widetilde{Z}$ and the linearization of $T$ at a balanced metric. For $H\in\M$, let $\pp_H$ be be the space of H-self-
adjoint endomorphisms of $H^0(X, E)$ that are induced by elements of $H^0(X,End(E))$.

\begin{lem}\label{lemlin}
Let $H\in\M$, let $a$ be an $H$-self-adjoint endomorphism of $H^0(X,E)$, and put $H_t(s,s')=H(e^{ta}s,s')$. Then the following hold.
\begin{enumerate}
\item $\frac{d}{dt}\big|_{t=0}\widetilde{Z}(H_t)=-\frac{rV}{N}\tr\big(a(A_H-I)\big)$. In particular, $H$ is a critical point of $\widetilde{Z}$ if and only if $H$ is balanced.
\item If $H$ is balanced, then
$$d_HT(a)=\frac{N}{Vr}\int_XP_x\,a\,P_x\,\frac{\omega^m}{m!},$$
where we identify tangent vectors at $H$ with $H$-self-adjoint endomorphisms.
\item If $H$ is balanced, then $d_HT$ is self-adjoint with respect to $\tr(ab)$, its eigenvalues lie in $(0,1]$, and the eigenvalue $1$ occurs exactly on $\pp_H$. Moreover, $\pp_H$ is the tangent space $T_H\B$ of $\B$ at $H$.
\end{enumerate}
\end{lem}

\begin{proof}
Let $\sigma_v$, $v\in E_x$, be as in the proof of Lemma \ref{partialorder}. By the proof of Lemma \ref{partialorder}, $|v|^2_{\FS(H_t)}$ is the minimum of $H_t(s,s)$ over the sections with $s(x)=v$, and for $t=0$ the minimum is attained at $\sigma_v$. Hence $|v|^2_{\FS(H_t)}\leq H_t(\sigma_v,\sigma_v)$ for all $t$, with equality at $t=0$. Since both sides are smooth in $t$, their derivatives agree at $t=0$, i.e.
$$\frac{d}{dt}\Big|_{t=0}|v|^2_{\FS(H_t)}=\frac{d}{dt}\Big|_{t=0}H_t(\sigma_v,\sigma_v)=H(a\sigma_v,\sigma_v).$$
By polarization, $\frac{d}{dt}\big|_{t=0}\FS(H_t)(v,w)=H(a\sigma_v,\sigma_w)=H(P_xaP_x\sigma_v,\sigma_w)$. Since $v\mapsto\sigma_v$ is an isometry from $(E_x,\FS(H))$ onto the range of $P_x$, we get
$$\frac{d}{dt}\Big|_{t=0}\log\det\FS(H_t)(x)=\tr(P_xaP_x)=\tr(aP_x).$$
Hence \eqref{eq1} and Lemma \ref{partialorder} give
$$\frac{d}{dt}\Big|_{t=0}Z(H_t)=-\int_X\tr(aP_x)\,\frac{\omega^m}{m!}=-\frac{Vr}{N}\tr(aA_H).$$
On the other hand, $\frac{d}{dt}\log\det H_t=\tr(a)$. This proves the formula in $(1)$. Since $A_H-I$ is $H$-self-adjoint, $H$ is a critical point of $\widetilde{Z}$ if and only if $A_H=I$, i.e. $T(H)=H$.

For $(2)$, we have
\begin{align*}
\frac{d}{dt}\Big|_{t=0}T(H_t)(s,s')&=\frac{N}{Vr}\int_XH\big(a\sigma_{s(x)},\sigma_{s'(x)}\big)\frac{\omega^m}{m!}\\
&=\frac{N}{Vr}\int_XH(aP_xs,P_xs')\,\frac{\omega^m}{m!},
\end{align*}
since $\sigma_{s(x)}=P_xs$. Since $T(H)=H$, the left hand side is $H(d_HT(a)s,s')$. This proves $(2)$.

For $(3)$, note that $\tr\big(b\,d_HT(a)\big)=\frac{N}{Vr}\int_X\tr(bP_xaP_x)\frac{\omega^m}{m!}$ is symmetric in $a$ and $b$. Moreover, since $A_H=I$, we have
\begin{align*}
\tr\big(a\,d_HT(a)\big)&=\frac{N}{Vr}\int_X\norm{P_xaP_x}^2\frac{\omega^m}{m!}\\
&\leq\frac{N}{Vr}\int_X\norm{P_xa}^2\frac{\omega^m}{m!}=\tr(a^2A_H)=\tr(a^2),
\end{align*}
with equality if and only if $P_xa(I-P_x)=0$ for all $x\in X$, i.e. $a$ commutes with every $P_x$. Hence the eigenvalues of $d_HT$ are at most $1$, and the eigenvalue $1$ occurs exactly on the space of those $a$ which preserve $K_x$ for every $x\in X$. As in the proof of Proposition \ref{nsprop1}, such an $a$ descends to a holomorphic endomorphism of $E$, i.e. $a\in\pp_H$. Conversely, every element of $\pp_H$ preserves $K_x$ for every $x \in X.$
\end{proof}

 Following Thompson, define 
\begin{equation}\label{thompsondistance}\dist_{T}(H,H^{'})=\inf\set{\lambda\geq 0 \mid e^{-\lambda}H\leq H'\leq
e^{\lambda}H}.\end{equation}
The following lemma is a simple observation that plays a crucial role in our argument. 
\begin{lem}\label{distlem}

For $H,H'\in \M$, we have
$$\dist(H,H')=\dist_T(H,H')=\inf\set{\lambda\geq 0 \mid e^{-\lambda}H\leq H'\leq
e^{\lambda}H}.$$
Moreover, 
\begin{enumerate}
\item $\dist$ is a metric on $\M$.
\item $\dist(e^{\alpha}H,e^{\alpha}H')=\dist(H,H')$
and $\dist(e^{\alpha}H,H)=|\alpha|$ for any $\alpha\in\R$.
\item For any $g\in GL(H^{0}(X,E))$, $\dist(g^{*}H,g^{*}H')=\dist(H,H')$, where
$g^{*}H(s,t)=H(gs,gt)$.
\item $H_{n}\rightarrow H$ in $\M$ if and only if
$\dist(H_{n},H)\rightarrow 0$.
\end{enumerate}

\end{lem}

\begin{proof}

 Let $s_{1},\dots,s_{N}$ be an $H$-orthonormal basis such that $H'(s_{i},s_{j})=\mu_{i}\delta_{ij}$. For $s=\sum_{i}c_{i}s_{i}$, we have
$$H(s,s)=\sum_{i}|c_{i}|^{2},\,\,\,\,\,\,\,\, H'(s,s)=\sum_{i}\mu_{i}|c_{i}|^{2}.$$
Therefore $e^{-\lambda}H\leq H'\leq e^{\lambda}H$ if and only if $e^{-\lambda}\leq \mu_{i}\leq e^{\lambda}$ for all $i$, i.e. $\lambda\geq \max_{i}|\log\mu_{i}|$. Hence $\dist_{T}(H,H')=\dist(H,H')$.

The eigenvalues of $H$ with respect to $H'$ are $\mu_{1}^{-1},...,
\mu_{N}^{-1}$. Hence $\dist$ is symmetric by Definition \ref{distance}, and $\dist(H,H')=0$
if and only if $\mu_{1}=...=\mu_{N}=1$, i.e. $H=H'$. If
$e^{-\lambda}H\leq H'\leq e^{\lambda}H$ and $e^{-\lambda'}H'\leq
H''\leq e^{\lambda'}H'$, then $e^{-\lambda-\lambda'}H\leq H''\leq
e^{\lambda+\lambda'}H$. This proves the triangle inequality and $(1)$.
Part $(2)$ is clear from Definition \ref{distance}, and $(3)$ holds since
$H\leq H'$ if and only if $g^{*}H\leq g^{*}H'$. For $(4)$, let
$s_{1},...,s_{N}$ be an $H$-orthonormal basis and let
$\mu_{1}^{(n)},...,\mu_{N}^{(n)}$ be the eigenvalues of
$[H_{n}(s_{i},s_{j})]$. Then $[H(s_{i},s_{j})]=I$ and
$$\big|[H_{n}(s_{i},s_{j})]-[H(s_{i},s_{j})]\big|_{op}=\max_{i}|\mu_{i}^{(n)}-1|,$$
which tends to zero if and only if $\dist(H_{n},H)=\max_{i}|\log
\mu_{i}^{(n)}|$ tends to zero.
\end{proof}

 Now we prove Proposition \ref{nsprop2}
 and Theorem \ref{mainthm}.

 \begin{proof}[Proof of Proposition \ref{nsprop2}]
     Let $H, H^{'} \in \M$ and let $\lambda=\dist(H,H^{'}).$ Then Lemma \ref{distlem} implies that $e^{-\lambda}H \leq H^{'} \leq e^{\lambda}H.$ By Lemma \ref{partialorder}, we also have
     $$e^{-\lambda}T(H)=T(e^{-\lambda}H)\leq T(H^{'})\leq T(e^{\lambda}H)=e^{\lambda}T(H).$$
     Applying Lemma \ref{distlem}, we get $$\dist(T(H),T(H^{'}))=\dist_T(T(H),T(H^{'}))\leq \dist(H,H^{'}).$$
 \end{proof}

\begin{proof}[Proof of Theorem \ref{mainthm}]
Let $H_{n}=T^{n}(H_{0}).$
We have 
\[\dist(H_{n+1},H_{\infty}) \leq \dist(H_{n},H_{\infty}) \leq \dots \leq \dist(H_{0},H_{\infty}).  \]
Therefore, $\{H_n\}$ is a bounded sequence. Therefore $\{H_{n}\}$ has a subsequence $H_{n_{j}}$ converging to some $\widetilde{H}_{\infty}\in \M$.

By Lemma \ref{firstlemma}, the sequence $\log\det H_{n}$ is non-increasing, and it is bounded from below. Hence $\log\det H_{n}$ converges to a real number $\ell$.  Since $T$ and $\log \det$ are continuous, we have
\[\log\det T(H_{\infty})=\lim_{j\to \infty} \log \det T(H_{n_j})=\lim_{j\to \infty} \log \det H_{n_j+1}=\log\det H_{\infty}.\]

Therefore,   $T(\widetilde{H}_{\infty})=\widetilde{H}_{\infty}$ by the equality case in Lemma \ref{firstlemma} and hence, $\widetilde{H}_{\infty}$ is a balanced metric. Proposition \ref{nsprop2} now implies that
$$\dist(H_{n+1},\widetilde{H}_{\infty})=\dist(T(H_{n}),T(\widetilde{H}_{\infty}))\leq \dist(H_{n},\widetilde{H}_{\infty}).$$
Hence $\dist(H_{n},\widetilde{H}_{\infty})$ is non-increasing. Since a subsequence of it tends to zero, the whole sequence tends to zero, i.e. $H_n\rightarrow \widetilde{H}_\infty$. Finally, by Proposition \ref{nsprop1}, there exists $\varphi\in\Aut(E)$, self-adjoint and positive with respect to $H_\infty$, such that $\widetilde{H}_\infty=\varphi^*H_\infty$.
\end{proof}

\begin{rem}\label{genrem}
The proof of convergence uses only two properties of the maps $\Hilb$ and $\FS$: the map $\Hilb$ is integration against a fixed positive measure, and $\FS(H)$ is the quotient metric induced by $H$ (Lemma \ref{partialorder}). Therefore the convergence of $T^{n}(H_{0})$ to a balanced metric holds if
the volume form $\frac{\omega^{m}}{m!}$ replaced by any finite positive measure $\nu$ on $X$ such that no nonzero section of $E$ vanishes on the support of $\nu$, with $V=\nu(X)$. The latter condition guarantees that $\Hilb(h)$ is positive definite for any $h\in\K$. Since $E$ is globally generated, $\tr P_{x}=r$ for every $x\in X$, and the trace identity in Lemma \ref{firstlemma} holds with the same proof. Hence, if a $\nu$-balanced metric exists, then $T_{\nu}^{n}(H_{0})$ converges to a $\nu$-balanced metric. Note that Proposition \ref{nsprop1} does not extend to this setting, and the limit need not lie in the $\Aut(E)$-orbit of a given $\nu$-balanced metric.

\end{rem}
\begin{rem}
In the case of line bundles, $T_{\nu}$ is the map defined by Donaldson in \cite{D3}. In \cite[Proposition 4]{D3}, Donaldson proves that $T_{\nu}^{n}(H_{0})$ converges under assumptions on $\nu$ which imply that the $\nu$-balanced metric is unique up to scaling. Our argument only uses the existence of a $\nu$-balanced metric. For example, let $\nu$ be the measure on $\CP^{2}$ with mass $\frac{2}{9}$ at each of three points of a line $\ell$ and mass $\frac{1}{3}$ at a point $p\notin\ell$. Then the $\nu$-balanced metrics on $\mathcal{O}(1)$ are not unique up to scaling, and $T_{\nu}^{n}(H_{0})$ still converges to one of them. For instance, let $\ell=\set{z_2=0}$, let the three points be $[1:\zeta^j:0]$, $j=0,1,2$, where $\zeta=e^{2\pi\sqrt{-1}/3}$, and let $p=[0:0:1]$. Then every inner product on $H^0(\CP^2,\mathcal{O}(1))$ for which $z_0,z_1,z_2$ are orthogonal and $\norm{z_0}=\norm{z_1}$ is $\nu$-balanced. In general, if $H$ is $\nu$-balanced, then $K_p^\perp$ is $H$-orthogonal to $K_q^\perp$ for $q\in\ell$, and rescaling $H$ on $K_p^\perp$ gives another $\nu$-balanced metric. Note that $\nu(\ell)/2=\nu(p)=\nu(\CP^2)/3$, so that $\nu$ violates the strict inequality in condition (2) of \cite[Proposition 4]{D3} only by equality. This case is relevant in practice, since in the numerical implementations of \cite{D3} and \cite{DKLR} the integral is replaced by a sum over finitely many points.
\end{rem}

\section{Rate of convergence and examples}\label{secrate}

In this section, we show that the convergence in Theorem \ref{mainthm} is exponential, and we study its rate. In Example \ref{excp}, we compute the rate for $\mathcal{O}(d)\otimes\C^{k}$ over $\CP^{m}$. In \cite[Theorem 1.2(1)]{IP}, Ioos and Polterovich show that for a stable bundle $E$, the rate for $E\otimes L^{p}$ is determined asymptotically by the first positive eigenvalue of the Laplacian on $\textrm{End}\,E$. In Proposition \ref{proppoly}, we extend their result to polystable bundles.

Using the argument of \cite[Section 4.2]{Se3}, we can prove a stronger version of Theorem \ref{mainthm}. We only give a sketch of proof. Suppose that $\B\neq\M$. Let $\dist_R$ be the Riemannian distance on $\M\cong GL(N,\C)/U(N)$, i.e. $\dist_R(H,H')=\big(\sum_i(\log\mu_i)^2\big)^{1/2}$ with $\mu_1,\dots,\mu_N$ as in Definition \ref{distance}. For $H\in\B$, let $\rho_*$ be the largest eigenvalue of $d_HT$ on the orthogonal complement of $T_H\B$. By Lemma \ref{lemlin}, $\rho_*\in(0,1)$, and by Lemma \ref{nslem4} and Proposition \ref{nsprop1} it does not depend on $H\in\B$. Then for any $R>0$, there is a constant $C=C(R)$ such that
\begin{equation}\label{optimalrate}\dist_R\big(T^n(H_0),\widetilde{H}_\infty\big)\leq C\rho_*^n\qquad\textrm{for all } n\geq0,\end{equation}
whenever $\dist_R(H_0,\B)\leq R$, where $\widetilde{H}_\infty=\lim_{n\to\infty}T^n(H_0)$. Moreover, $\rho_*$ cannot be replaced by any smaller number. Indeed, the proofs of \cite[Theorem 1.1 and Proposition 4.7]{Se3} only use the following facts, which hold in the present setting. The map $T$ is smooth, since by Lemma \ref{partialorder} it is given by an integral of the projections $P_x$, which depend smoothly on $H$. The group $\Aut(E)$ acts on $\M$ by isometries and transitively on $\B$, and $T$ is equivariant (Lemma \ref{nslem4} and Proposition \ref{nsprop1}). By Proposition \ref{nsprop1} and the arguments of \cite[Corollary 2.9 and Proposition 2.11]{Se3}, the set $\B=\set{(e^{a})^*H_\infty\mid a\in\pp_{H_\infty}}$ is a closed totally geodesic submanifold of $\M$. The functional $\widetilde{Z}$ is convex, non-increasing along the orbits of $T$ and $\Aut(E)$-invariant (Theorem \ref{thirdthm}, Lemma \ref{firstlemma} and Proposition \ref{propinv}). By Lemma \ref{lemlin} and the argument of \cite[Theorem 4.2(1)]{Se3}, we have $d_HT=I-\frac{N}{rV}\nabla^2\widetilde{Z}(H)$ for $H\in\B$. Finally, $T^n(H_0)$ converges to a point of $\B$ by Theorem \ref{mainthm}.

\begin{ex}\label{excp}
Let $X=\CP^m$ with the Fubini-Study form $\omega_{FS}$, and let $E=\mathcal{O}(d)\otimes\C^k$ with $d\geq1$ and $k\geq1$. Put $W=H^0(\CP^m,\mathcal{O}(d))$ and $N_L=\dim W$, so that $H^0(X,E)=W\otimes\C^k$, $N=kN_L$ and $r=k$. Let $H_L$ be a $U(m+1)$-invariant inner product on $W$, and let $H=H_L\otimes\langle\cdot,\cdot\rangle_{\C^k}$. For $x\in X$, the subspace $K_x$ for $E$ is $K^L_x\otimes\C^k$, where $K^L_x\subseteq W$ is the corresponding subspace for $\mathcal{O}(d)$. Hence $P_x=P^L_x\otimes I$, where $P^L_x$ is the $H_L$-orthogonal projection onto the line $(K^L_x)^\perp$. Therefore $A_H=A^L\otimes I$, where $A^L=\frac{N_L}{V}\int_XP^L_x\frac{\omega_{FS}^m}{m!}$. Since $H_L$ and $\omega_{FS}$ are $U(m+1)$-invariant, $A^L$ commutes with the action of $U(m+1)$ on $W$. Since $W$ is an irreducible representation of $U(m+1)$ and $\tr A^L=N_L$, we get $A^L=I$. Hence $H$ is balanced.

For $b\in\Herm_{H_L}(W)$, put $\sigma_L(b)(x)=\tr(P^L_xb)$, and for a real function $f$ on $X$ define $Q_L(f)\in\Herm_{H_L}(W)$ by
$$H_L\big(Q_L(f)s,s'\big)=\frac{N_L}{V}\int_Xf\,\FS(H_L)(s,s')\,\frac{\omega_{FS}^m}{m!}.$$
Then $\tr\big(Q_L(f)b\big)=\frac{N_L}{V}\int_Xf\,\sigma_L(b)\,\frac{\omega_{FS}^m}{m!}$. Since $P^L_x$ has rank one, we have $P_x(b\otimes c)P_x=\sigma_L(b)(x)\,P^L_x\otimes c$ for $c\in\Herm(\C^k)$. Hence, by Lemma \ref{lemlin}(2),
\begin{align*}
\big\langle d_HT(b\otimes c),b'\otimes c'\big\rangle_H&=\frac{N_L}{V}\int_X\sigma_L(b)\,\sigma_L(b')\,\frac{\omega_{FS}^m}{m!}\cdot\tr(cc')\\
&=\tr\big(Q_L(\sigma_L(b))\,b'\big)\,\tr(cc').
\end{align*}
Therefore $d_HT=(Q_L\circ\sigma_L)\otimes\textrm{Id}$ on $\Herm_H=\Herm_{H_L}(W)\otimes_\R\Herm(\C^k)$. Note that $\sigma_L$ and $Q_L$ are the operators $\sigma$ and $Q$ of \cite[Remark 4.3]{Se3} for $(\CP^m,\mathcal{O}(d))$ at the balanced inner product $H_L$, since the volume forms $\frac{\omega_{FS}^m}{m!}$ and $\frac{(d\,\omega_{FS})^m}{m!}$ are proportional. Let $U_\ell\subseteq\Herm_{H_L}(W)$, $0\leq\ell\leq d$, be the subspaces of \cite[Section 5.1]{Se3}, so that $\Herm_{H_L}(W)=\bigoplus_{\ell=0}^dU_\ell$. By \cite[Lemma 5.2]{Se3} ( see also \cite[Lemma 3.7]{Se2} )and the proof of \cite[Lemma 5.1]{Se3}, $Q_L\circ\sigma_L$ acts on $U_\ell$ by the scalar
$$\beta_\ell=\frac{d!\,(d+m)!}{(d-\ell)!\,(d+m+\ell)!}.$$
Since $\beta_0=1$ and $\beta_{\ell+1}/\beta_\ell=(d-\ell)/(d+m+\ell+1)<1$, Lemma \ref{lemlin} gives $T_H\B=\pp_H=U_0\otimes\Herm(\C^k)$ and
$$\rho_*=\beta_1=\frac{d}{d+m+1}.$$
\end{ex}

\begin{rem}\label{remcompare}
In Example \ref{excp}, we have $1-\rho_*=\frac{m+1}{d+m+1}$. For Donaldson's map $T$ on $(\CP^m,\mathcal{O}(d))$, the corresponding number satisfies $1-\rho_*\sim\frac{(m+2)(m+3)}{d^2}$ as $d\to\infty$ (\cite{Se3}). The difference comes from the fact that the volume form is fixed here. Indeed, the directions $U_1\otimes\Herm(\C^k)$ are normal to $\B$, while for Donaldson's map the directions $U_1$ are tangent to the set of balanced metrics. Moreover, for a K\"ahler form in $c_1(\mathcal{O}(1))$, the first positive eigenvalue of the Laplacian $2\bar\partial^*\bar\partial$ on $\textrm{End}(\mathcal{O}^{\oplus k})$ is $4\pi(m+1)$. Hence Example \ref{excp} agrees with the asymptotic formula $1-\frac{\lambda_1^E}{4\pi p}+o(p^{-1})$ of \cite[Theorem 1.2(1)]{IP} for $E=\mathcal{O}^{\oplus k}$ and $p=d$, where $\lambda_1^E$ is the first positive eigenvalue of $2\bar\partial^*\bar\partial$ on $\textrm{End}\,E$ (\cite[(1.12)]{IP}), although $\mathcal{O}^{\oplus k}$ is not stable for $k\geq2$. Note that in \cite[(1.13)]{IP} the eigenvalues are listed with multiplicity, so that the kernel $\mathfrak{gl}(k)$ contributes $k^2$ zero eigenvalues. In Proposition \ref{proppoly} below, we show that this asymptotic formula holds for every polystable bundle $E$ such that $E\otimes L^{p}$ admits balanced metrics for $p\gg0$.
\end{rem}

Throughout the rest of this section, let $L$ be an ample line bundle on $X$ such that $\omega\in c_{1}(L)$, and let $h_{L}$ be a Hermitian metric on $L$ whose curvature form is $-2\pi\sqrt{-1}\,\omega$. For $p\geq0$, put $E(p)=E\otimes L^{p}$. For $p\gg0$, the bundle $E(p)$ is very ample, and we apply the results above to $E(p)$ with the fixed volume form $\frac{\omega^{m}}{m!}$. We denote by $T_{p}$, $\M_{p}$ and $\B_{p}$ the corresponding objects, we put $N_{p}=\dim H^{0}(X,E(p))$. Assume that $E(p)$ admits balanced metrics. We denote by $\rho_*(p)$ the largest eigenvalue of $d_HT_{p}$ on the orthogonal complement of $T_H\B_{p}$ for $H \in \B_p.$ Note that  $\rho_*(p)$ does not depend on $H$ since $T_p$ is equivariant under the action  of $\Aut(E).$ By  \eqref{optimalrate},  $\rho_*(p)$ is the optimal rate of convergence for $E(p)$ by $\rho_*(p)$. Since $\textrm{End}\,E(p)=\textrm{End}\,E$, the group $\Aut(E)$ and algebra $\A=H^{0}(X,\textrm{End}\,E)$ act on $H^{0}(X,E(p))$ for all $p$.

For a Hermitian metric $h$ on $E$, let $\square_{h}=2\bar\partial^{*}\bar\partial$ be the Laplacian acting on $C^{\infty}(X,\textrm{End}\,E)$, where $\bar\partial^{*}$ is the formal adjoint of $\bar\partial$ with respect to the $L^{2}(h,\frac{\omega^{m}}{m!})$ (\cite[(1.12)]{IP}). It is a non-negative self-adjoint elliptic operator whose kernel is $\ker\bar\partial=\A$. 
Note that for $\varphi\in\Aut(E)$, the operators $\square_{\varphi^{*}h}$ and $\square_{h}$ are conjugate and therefore, they have the same eigenvalues. The proof of the following 
is very similar to the case of stable bundles treated in \cite{IP}. 

\begin{prop}\label{proppoly}
Let $E=F_{1}\oplus\dots\oplus F_{l}$, where $F_{1},\dots,F_{l}$ are Mumford stable with respect to $L$. Suppose that
\begin{equation}\label{eqhilb}
\frac{\chi(X,F_{i}\otimes L^{p})}{\textrm{rk}(F_{i})}=\frac{\chi(X,E\otimes L^{p})}{\textrm{rk}(E)}
\end{equation}
for $1\leq i\leq l$ and $p \gg 0$. Then $E(p)$ admits a balanced metric for $p\gg0$, and
$$\rho_*(p)=1-\frac{\lambda_{1}^{E}}{4\pi p}+o(p^{-1})\qquad\textrm{as } p\rightarrow\infty,$$
where $\lambda_{1}^{E}$ is the first positive eigenvalue of $\square_{h_{E}}$, $h_E$ is a Hermitian-Einstein metric on $E$.
\end{prop}

Note that $\lambda_{1}^{E}$ does not depend on the choice of the Hermitian-Einstein metric on $E$.  For Mumford stable bundle (i.e. $l=1$)  Proposition \ref{proppoly} is the estimate (1.22) of \cite[Theorem 1.2(1)]{IP}. For $E=\mathcal{O}^{\oplus k}$ over $\CP^{m}$ and $L=\mathcal{O}(1)$, it agrees with Example \ref{excp} (see Remark \ref{remcompare}).

\begin{proof}
Put $r=\textrm{rk}(E)$, $r_{i}=\textrm{rk}(F_{i})$ and $N_{i,p}=\dim H^{0}(X,F_{i}(p))$. Since $L$ is ample, \eqref{eqhilb} gives $N_{i,p}/r_{i}=N_{p}/r$ for $p\gg0$.

 By the theorem of Wang (\cite[Theorem 1.2]{W2}, see also \cite[Theorem 4.5]{IP}), for $p\gg0$ there is a balanced metric $H_{i,p}$ on $H^{0}(X,F_{i}(p))$ such that 
$$h_p=h_{L}^{-p}\otimes\FS(H_{p})\To e^{f}h_{E}=:h_{\infty}\qquad\textrm{in } C^{\infty} \textrm{ as } p\rightarrow\infty,$$
where $f$ is a smooth function on $X$. Here,
$$H_{p}=H_{1,p}\oplus\dots\oplus H_{l,p}\quad\textrm{on}\quad H^{0}(X,E(p))=\bigoplus_{i=1}^{l}H^{0}(X,F_{i}(p)).$$

It is easy to prove that the metric $H_p$ is balanced. 
By \eqref{optimalrate}, we may compute $\rho_*(p)$ at $H_{p}$.

Let $d=\dim_{\C}\A$.
Therefore $\pp_{H_{p}}=\set{a\in\A\mid a^{*}=a}$ has real dimension $d$ which implies that eigenvalue $1$ of $d_{H_p}T_p$ has multiplicity $d.$ Here $a^*$ is the $H_p$-adjoint of $a.$  Let
$$1=\beta_{0,p}\geq\beta_{1,p}\geq\beta_{2,p}\geq\dots>0$$
be the eigenvalues of $d_{H_{p}}T_{p}$ on the space of $H_p$-self-adjoint endomorphisms, counted with multiplicity. By Lemma \ref{lemlin}(3), we get
\begin{equation}\label{eqbeta}
\beta_{0,p}=\dots=\beta_{d-1,p}=1>\beta_{d,p}=\rho_*(p).
\end{equation}

By \cite[Corollary 4.4]{IP}, we have $d_{H_{p}}T_{p}=T_{p}\circ T_{p}^{*}$, where $T_{p}$ is the Berezin-Toeplitz quantization map associated with the volume form $\frac{\omega^{m}}{m!}$ and the metric $\FS(H_{p})$ on $E(p)$. Note that $\FS(H_{p})=h_{L}^{p}\otimes h_{p}$, and that the identity $d_{H_{p}}T_{p}=T_{p}\circ T_{p}^{*}$ agrees with Lemma \ref{lemlin}(2). Hence $\beta_{k,p}$, $k\geq0$, are the positive eigenvalues of the Berezin transform $T_{p}^{*}\circ T_{p}$ of $(E,h_{p})$ over $(X,p\,\omega)$ (\cite[Definition 3.1]{IP}), counted with multiplicity. Let
$$0=\lambda_{0}\leq\lambda_{1}\leq\lambda_{2}\leq\dots$$
be the eigenvalues of $\square_{h_{\infty}}=\square_{h_{E}}$, counted with multiplicity. We claim that for every fixed $k\geq0$,
\begin{equation}\label{eqasym}
\beta_{k,p}=1-\frac{\lambda_{k}}{4\pi p}+o(p^{-1})\qquad\textrm{as } p\rightarrow\infty.
\end{equation}
If the metrics $h_{p}$ did not depend on $p$, then \eqref{eqasym} would be \cite[Theorem 1.1]{IP}. In our case, the metrics $h_{p}$ converge to $h_{\infty}$ in $C^{\infty}$, and the constants in \cite[Propositions 3.4 and 3.6]{IP} are uniform in the $C^{l}$-norms of the metric on $E$ for some $l$. Hence the argument of the proof of \cite[Theorem 4.17]{IP}, which is also the argument of the proof of \cite[Theorem 1.2(1)]{IP}, gives \eqref{eqasym}. Note that this part of \cite{IP} does not use the stability of $E$. 
Finally, the kernel of $\square_{h_{\infty}}$ is $\A$. Hence $\lambda_{0}=\dots=\lambda_{d-1}=0<\lambda_{d}=\lambda_{1}^{E}$. Applying \eqref{eqasym} with $k=d$ and using \eqref{eqbeta}, we get
$$\rho_*(p)=\beta_{d,p}=1-\frac{\lambda_{1}^{E}}{4\pi p}+o(p^{-1}).$$
\end{proof}

\begin{rem}\label{rempoly}
\begin{enumerate}
\item Let $E$ be polystable, i.e. $E=F_{1}\oplus\dots\oplus F_{l}$, where the bundles $F_{i}$ are stable and have the same slope. If $E(p)$ admits a balanced metric for infinitely many $p$, then Corollary \ref{cordirect} gives $N_{i,p}/r_{i}=N_{p}/r$ for infinitely many $p$, and therefore \eqref{eqhilb} holds for all $p$. Hence Proposition \ref{proppoly} applies to every polystable bundle $E$ such that $E(p)$ admits balanced metrics for $p\gg0$, i.e. to every polystable bundle for which $\rho_*(p)$ is defined for $p\gg0$.
\item The condition that the slopes agree does not suffice. For example, let $X=\CP^{1}\times\CP^{1}$, $L=\mathcal{O}(1,1)$ and $E=\mathcal{O}(2,0)\oplus\mathcal{O}(1,1)$. Then $E$ is polystable, but
$$\chi(X,\mathcal{O}(p+2,p))=(p+3)(p+1)\neq(p+2)^{2}=\chi(X,\mathcal{O}(p+1,p+1)).$$
Hence Corollary \ref{cordirect} gives $F(\beta_{1})=\frac{1}{2}$ for $E(p)$ and every $p\geq1$. Therefore $E(p)$ admits no balanced metric, and its Gieseker point is unstable (Remark \ref{remgieseker}). In accordance with Remark \ref{remslope}, the leading coefficient of $F(\beta_1)$ vanishes, but $F(\beta_1)$ does not.
\item The proof shows that at any balanced metric $H$, the algebra $\A$ is closed under taking $H$-adjoints, and the multiplicity of the eigenvalue $1$ of $d_{H}T$ is $\dim_{\C}\A$. For simple bundles, the relation between the linearization at balanced metrics and the Laplacian on $\textrm{End}\,E$ is also studied in \cite{KMS} (see \cite[Remark 4.20]{IP}).
\end{enumerate}
\end{rem}

In the stable case, $\rho_*(p)$ is controlled by the eigenvalues of the Laplacian on $\textrm{End}\,E$. If $E$ is polystable, then $\textrm{End}\,E$ contains the bundles $\textrm{Hom}(F_{j},F_{i})$, $i\neq j$, and the following example shows that they can determine the rate.

\begin{ex}\label{exP1P1}
Let $X=\CP^{1}\times\CP^{1}$, $L=\mathcal{O}(1,1)$ and $\omega=\pi_{1}^{*}\omega_{FS}+\pi_{2}^{*}\omega_{FS}$, where $\omega_{FS}\in c_{1}(\mathcal{O}(1))$ is the Fubini-Study form, so that $V=1$. For $k\geq1$, let $E=\mathcal{O}(k,0)\oplus\mathcal{O}(0,k)$ and hence, we have $\A=\C\oplus\C$. We compute $\rho_*(p)$ for $p\geq1$.

Let $W_{j}=H^{0}(\CP^{1},\mathcal{O}(j))$ and $n_{j}=j+1$. Let $H_{W_{j}}$ be a $U(2)$-invariant inner product on $W_{j}$, and let $P^{(j)}_{z}$, $z\in\CP^{1}$, be the corresponding projections. By Example \ref{excp} with $m=1$, $H_{W_j}$ is balanced, i.e.
$$n_{j}\int_{\CP^{1}}P^{(j)}_{z}\,\omega_{FS}=I.$$
We have
$$H^{0}(X,E(p))=(W_{p+k}\otimes W_{p})\oplus(W_{p}\otimes W_{p+k}),$$
and we put $H_{p}=(H_{W_{p+k}}\otimes H_{W_{p}})\oplus(H_{W_{p}}\otimes H_{W_{p+k}})$. Then
$$P_{(z,w)}=(P^{(p+k)}_{z}\otimes P^{(p)}_{w})\oplus(P^{(p)}_{z}\otimes P^{(p+k)}_{w}),$$
and therefore $H_p$ is balanced. Since $\frac{\omega^{2}}{2}=\pi_{1}^{*}\omega_{FS}\wedge\pi_{2}^{*}\omega_{FS}$ and $\frac{N_{p}}{Vr}=n_{p+k}n_{p}$, Lemma \ref{lemlin}(2) shows that $d_{H_{p}}T_{p}$ preserves the decomposition of an endomorphism of $H^{0}(X,E(p))$ into four blocks, and that it acts by $Q_{p+k}\otimes Q_{p}$ and $Q_{p}\otimes Q_{p+k}$ on the diagonal blocks, and by $R\otimes R'$ and $R'\otimes R$ on
$$\textrm{Hom}(W_{p}\otimes W_{p+k},W_{p+k}\otimes W_{p})=\textrm{Hom}(W_{p},W_{p+k})\otimes\textrm{Hom}(W_{p+k},W_{p})$$
and on the other off-diagonal block. Here
\begin{align*}
Q_{j}(b)&=n_{j}\int_{\CP^{1}}P^{(j)}_{z}\,b\,P^{(j)}_{z}\,\omega_{FS},& b&\in\textrm{End}(W_{j}),\\
R(b)&=n_{p+k}\int_{\CP^{1}}P^{(p+k)}_{z}\,b\,P^{(p)}_{z}\,\omega_{FS},& b&\in\textrm{Hom}(W_{p},W_{p+k}),\\
R'(b)&=n_{p}\int_{\CP^{1}}P^{(p)}_{z}\,b\,P^{(p+k)}_{z}\,\omega_{FS},& b&\in\textrm{Hom}(W_{p+k},W_{p}).
\end{align*}
These operators are self-adjoint and non-negative with respect to the Hilbert-Schmidt inner products. By Example \ref{excp}, the eigenvalues of $Q_{j}$ are
$$\beta_{\ell}(j)=\frac{j!\,(j+1)!}{(j-\ell)!\,(j+\ell+1)!},\qquad 0\leq\ell\leq j,$$
with multiplicity $2\ell+1$. Hence the eigenvalue $1$ occurs on the diagonal blocks with multiplicity $2=\dim\A$, and the largest of the other eigenvalues on the diagonal blocks is $\beta_{1}(p+k)=\frac{p+k}{p+k+2}$. On the other hand, for $b\in\textrm{Hom}(W_{p},W_{p+k})$ we have
\begin{align*}
\tr\big(b^{*}R(b)\big)&=n_{p+k}\int_{\CP^{1}}\norm{P^{(p+k)}_{z}bP^{(p)}_{z}}^{2}\omega_{FS}\\
&\leq n_{p+k}\int_{\CP^{1}}\norm{P^{(p+k)}_{z}b}^{2}\omega_{FS}=\tr(b^{*}b),
\end{align*}
with equality if and only if $b$ maps the sections vanishing at $z$ to sections vanishing at $z$ for every $z$, i.e. $b$ is the multiplication by an element of $W_{k}$. Similarly, for $b\in\textrm{Hom}(W_{p+k},W_{p})$ we have
$$\tr\big(b^{*}R'(b)\big)\leq n_{p}\int_{\CP^{1}}\norm{bP^{(p+k)}_{z}}^{2}\omega_{FS}=\frac{n_{p}}{n_{p+k}}\tr(b^{*}b),$$
with equality if and only if $b^{*}$ is the multiplication by an element of $W_{k}$. Hence the largest eigenvalues of $R$ and $R'$ are $1$ and $\frac{p+1}{p+k+1}$, both with multiplicity $k+1$. Therefore the largest eigenvalue on the off-diagonal blocks is $\frac{p+1}{p+k+1}$, with real multiplicity $2(k+1)^{2}$ on self-adjoint endomorphisms, and
$$\rho_*(p)=\max\set{\frac{p+k}{p+k+2},\,\frac{p+1}{p+k+1}}.$$
In particular, $1-\rho_*(p)=\min\set{2,k}\,p^{-1}+O(p^{-2})$. This agrees with Proposition \ref{proppoly}. Indeed, $\textrm{End}\,E=\mathcal{O}\oplus\mathcal{O}\oplus\mathcal{O}(k,-k)\oplus\mathcal{O}(-k,k)$. The first positive eigenvalue of the Laplacian on functions on $(X,\omega)$ is $8\pi$. On $\mathcal{O}(k,-k)=\pi_{1}^{*}\mathcal{O}(k)\otimes\pi_{2}^{*}\mathcal{O}(-k)$, the operator $2\bar\partial^{*}\bar\partial$ is the sum of the corresponding operators along the two factors. Along the first factor, its kernel is $W_{k}$. Along the second factor, the Bochner-Kodaira-Nakano identity gives $2\bar\partial^{*}\bar\partial=2\partial^{*}\partial+4\pi k$ on the sections of $\mathcal{O}(-k)$, and the kernel of $\partial$ consists of the sections $s^{*}$, $s\in W_{k}$, where $s^{*}$ is the dual of $s$ with respect to the Hermitian metric. Hence the first eigenvalue of $2\bar\partial^{*}\bar\partial$ on $\mathcal{O}(k,-k)$ is $4\pi k$, with multiplicity $(k+1)^{2}$, and $\lambda_{1}^{E}=4\pi\min\set{2,k}$.

For $k=1$, we get $\rho_*(p)=\frac{p+1}{p+2}$, and the slowest directions are off-diagonal. On the other hand, for each summand $F_{i}(p)$ the number $\rho_*$ is $\beta_{1}(p+1)=\frac{p+1}{p+3}$. Hence for $E=\mathcal{O}(1,0)\oplus\mathcal{O}(0,1)$, the iteration converges more slowly than for its summands, and the rate is determined by the eigenvalues of the Laplacian on $\textrm{Hom}(F_{2},F_{1})$ and $\textrm{Hom}(F_{1},F_{2})$.
\end{ex}

\appendix
\section{Invariance of $\widetilde{Z}$ under $\Aut(E)$}\label{appinv}

In \cite[Lemma 2.7]{Se3}, we show that the functional $\widetilde{Z}$ of a polarized manifold is invariant under the automorphism group, provided that a balanced metric exists. In \cite[Appendix A]{Se3}, we show that this assumption cannot be removed, and we relate the invariance to the obstruction to asymptotic Chow semistability of Mabuchi and Futaki (\cite{Ma}, \cite{Fu2}). In this appendix, we discuss the analogue for vector bundles. In this case, the obstruction is given by an explicit formula, and the Chow point is replaced by the Gieseker point. Throughout this appendix, we do not assume that $E$ admits a balanced metric.

Let $\A=H^0(X,\textrm{End}\,E)$ be the Lie algebra of connected complex Lie group $\Aut(E)$. For $\beta\in\A$, the function $x\mapsto\tr(\beta_x)$ is holomorphic on $X$, hence constant. We denote this constant by $\tr_E(\beta)$. Similarly, for $\varphi\in\Aut(E)$, we denote the constant $\det(\varphi_x)$ by $\det\nolimits_E(\varphi)$. We denote by $\tr_{H^0}(\beta)$ and $\det\nolimits_{H^0}(\varphi)$ the trace and the determinant of the action on $H^0(X,E)$. Define the complex linear functional $F:\A\to\C$ by
$$F(\beta)=\frac{N}{r}\tr_E(\beta)-\tr_{H^0}(\beta).$$
Note that $F$ vanishes on the scalars and therefore, it is identically zero for asimple bundle $E$.

\begin{prop}\label{propinv}
For any $\varphi\in\Aut(E)$ and $H\in\M$, we have
$$\widetilde{Z}(\varphi^*H)-\widetilde{Z}(H)=2V\Big(\frac{r}{N}\log\big|\det\nolimits_{H^0}(\varphi)\big|-\log\big|\det\nolimits_E(\varphi)\big|\Big).$$
Moreover, the following hold.
\begin{enumerate}
\item For $\beta\in\A$ and $t\in\R$, we have
$$\widetilde{Z}\big((e^{t\beta})^*H\big)-\widetilde{Z}(H)=-\frac{2rV}{N}\,t\,\textrm{Re}\,F(\beta).$$
\item The functional $\widetilde{Z}$ is $\Aut(E)$-invariant if and only if $F=0$.
\item For any $H\in\M$, we have $F(\beta)=\tr\big(\beta(A_H-I)\big)$, where $\beta$ acts on $H^0(X,E)$ and $A_H$ is as in Lemma \ref{partialorder}. In particular, if $E$ admits a balanced metric, then $F=0$.
\end{enumerate}
\end{prop}

\begin{proof}
Let $h=\FS(H)$. By Lemma \ref{nslem4}, we have \[\det(\varphi^*h)=|\det\nolimits_E(\varphi)|^2\det(h). \] Integrating \eqref{eq1} along the path $h_t=e^{t\psi}h$, where $\psi$ is the endomorphism of $E$ with $\varphi^*h=e^{\psi}h$, we get
$$I(\varphi^*h)-I(h)=\int_X\log\big|\det\nolimits_E(\varphi)\big|^2\,\frac{\omega^m}{m!}=2V\log\big|\det\nolimits_E(\varphi)\big|.$$
This together with $\log\det(\varphi^*H)-\log\det(H)=2\log|\det\nolimits_{H^0}(\varphi)|$ proves the formula. For $(1)$, note that $\det\nolimits_E(e^{t\beta})=e^{t\tr_E(\beta)}$ and $\det\nolimits_{H^0}(e^{t\beta})=e^{t\tr_{H^0}(\beta)}$.

If $\widetilde{Z}$ is $\Aut(E)$-invariant, then $(1)$ implies that $F=0$, since $F$ is complex linear. Conversely, suppose that $F=0$. By the formula, $\varphi\mapsto\widetilde{Z}(\varphi^*H)-\widetilde{Z}(H)$ is a homomorphism from $\Aut(E)$ to $\R$, and by $(1)$ it vanishes on $e^{t\beta}$ for every $\beta\in\A$. Since $\Aut(E)$ is connected, it vanishes on $\Aut(E)$. This proves $(2)$.

For $(3)$, let $x\in X$, and let $\sigma_v$, $v\in E_x$, be as in the proof of Lemma \ref{partialorder}. Since $P_xs=\sigma_{s(x)}$ for every $s\in H^0(X,E)$, we have $P_x\beta\sigma_v=\sigma_{\beta(v)}$. Hence the restriction of $P_x\beta P_x$ to the range of $P_x$ is conjugate to $\beta_x$, and $\tr(\beta P_x)=\tr(P_x\beta P_x)=\tr_E(\beta)$. Therefore
$$\tr(\beta A_H)=\frac{N}{Vr}\int_X\tr(\beta P_x)\,\frac{\omega^m}{m!}=\frac{N}{r}\tr_E(\beta),$$
which gives $F(\beta)=\tr\big(\beta(A_H-I)\big)$. If $H$ is balanced, then $A_H=I$ by Lemma \ref{partialorder}, and therefore $F=0$.
\end{proof}

Note that $(1)$ and $(3)$ are the analogues of \cite[Proposition A.1]{Se3}. In particular, if $F\neq0$, then $E$ admits no balanced metric, and $\widetilde{Z}$ is not bounded from below on $\M$. Since $F$ vanishes in the case of simple bundle $E$,the obstruction is only relevant for non-simple bundles.

\begin{cor}\label{cordirect}
Suppose that $E=E_1\oplus\dots\oplus E_l$, and let $r_i$ and $N_i$ be the rank of $E_i$ and the dimension of $H^0(X,E_i)$. Let $\beta_i\in\A$ be the projection onto $E_i$. Then
$$F(\beta_i)=r_i\Big(\frac{N}{r}-\frac{N_i}{r_i}\Big).$$
In particular, if $E$ admits a balanced metric, then $N_i/r_i=N/r$ for all $i$.
\end{cor}

\begin{proof}
We have $\tr_E(\beta_i)=r_i$ and $\tr_{H^0}(\beta_i)=N_i$. Hence the claim follows from Proposition \ref{propinv}.
\end{proof}

\begin{rem}\label{remgieseker}
The functional $F$ is related to the Gieseker point $\T(E)$ in the same way as the obstruction of \cite[Appendix A]{Se3} is related to the Chow point. Indeed, for $\varphi\in\Aut(E)$ and $s_1,\dots,s_r\in H^0(X,E)$, we have
$$\T(E)(\varphi s_1\wedge\dots\wedge\varphi s_r)=\det\nolimits_E(\varphi)\,\T(E)(s_1\wedge\dots\wedge s_r).$$
Let $g=\lambda\varphi$, where $\lambda^N=\det\nolimits_{H^0}(\varphi)^{-1}$, so that $\det g=1$. Then $g$ preserves the line spanned by $\T(E)$, and it acts on this line by a scalar whose absolute value is $|\chi(\varphi)|$ or $|\chi(\varphi)|^{-1}$, depending on the convention for the action of $GL(N)$, where
$$|\chi(\varphi)|=\frac{|\det\nolimits_{H^0}(\varphi)|^{r/N}}{|\det\nolimits_E(\varphi)|}.$$
By Proposition \ref{propinv}, we have $\widetilde{Z}(\varphi^*H)-\widetilde{Z}(H)=2V\log|\chi(\varphi)|$. Suppose that $F\neq0$, and choose $\beta\in\A$ with $\textrm{Re}\,F(\beta)\neq0$. Then the scalar by which the element of $SL(N)$ associated with $e^{t\beta}$ acts on $\T(E)$ tends to $0$ as $t\to\infty$ or as $t\to-\infty$. Therefore $0$ lies in the closure of the $SL(N)$-orbit of $\T(E)$, i.e. the Gieseker point of $E$ is unstable. In particular, $F=0$ if the Gieseker point of $E$ is semistable. This is the analogue for vector bundles of the obstruction to Chow semistability of Mabuchi and Futaki (\cite{Ma}, \cite{Fu2}).
\end{rem}

\begin{rem}\label{remslope}
Let $L$ be an ample line bundle on $X$, and for $k\geq1$ let $F_k$ be the functional $F$ for $E\otimes L^k$. Note that $\A$ and $\Aut(E)$ do not change if $E$ is replaced by $E\otimes L^k$. Suppose that $E=E_1\oplus\dots\oplus E_l$. By Corollary \ref{cordirect} and the theorem of Riemann--Roch, for $k\gg0$ we have
$$F_k(\beta_i)=r_i\big(\mu(E)-\mu(E_i)\big)\frac{k^{m-1}}{(m-1)!}+O(k^{m-2}),$$
where $\mu(E')=\deg_L(E')/\textrm{rk}(E')$ is the slope with respect to $L$, and $\deg_L(E')=c_1(E')\cdot c_1(L)^{m-1}$. Hence the leading coefficient of $F_k$ is given by the differences of slopes. It vanishes if $E$ admits a Hermitian-Einstein metric with respect to a K\"ahler form in $c_1(L)$. Indeed, in this case $E$ is polystable (\cite{D1}, \cite{D2}, \cite{UY}), and every direct summand of a polystable bundle has the same slope as the bundle, since it is both a subbundle and a quotient bundle. This is the analogue of the fact that the first of the integral invariants of Futaki is the classical Futaki invariant, which obstructs the existence of K\"ahler-Einstein metrics (\cite{Fu1}, \cite{Fu2}, see also \cite[Remark A.2]{Se3}).
\end{rem}

\begin{ex}\label{exsplit}
Let $X=\CP^1$ and $E=\mathcal{O}(a)\oplus\mathcal{O}(b)$ with $a,b\geq1$. Then $N=a+b+2$ and $r=2$, and Corollary \ref{cordirect} gives
$$F(\beta_1)=\frac{b-a}{2}.$$
Explicitly, for $\varphi=\lambda\,\textrm{Id}_{\mathcal{O}(a)}\oplus\textrm{Id}_{\mathcal{O}(b)}$ with $\lambda\in\C^*$, Proposition \ref{propinv} gives
$$\widetilde{Z}(\varphi^*H)-\widetilde{Z}(H)=\frac{2V(a-b)}{a+b+2}\log|\lambda|.$$
Hence, if $a\neq b$, the functional $\widetilde{Z}$ is neither $\Aut(E)$-invariant nor bounded from below, $E$ admits no balanced metric, and its Gieseker point is unstable. The same holds for all twists $E\otimes\mathcal{O}(j)$, $j\geq0$, since $F(\beta_1)$ does not depend on $j$, in accordance with Remark \ref{remslope}. On the other hand, for $E=\mathcal{O}(d)\otimes\C^k$ on $\CP^m$ as in Example \ref{excp}, we have $\A=\textrm{End}(\C^k)$, $\tr_{H^0}(\beta)=N_L\tr(\beta)$ and $N/r=N_L$. Hence $F=0$, in accordance with the existence of balanced metrics.
\end{ex}


\begin{thebibliography}{HHHH}

\bibitem[BCCT]{BCCT} J. Bennett, A. Carbery, M. Christ and T. Tao, The Brascamp-Lieb inequalities: finiteness, structure and extremals, Geom.\ Funct.\ Anal. \textbf{17} (2008), 1343--1415.



\bibitem[D1]{D1} S. K. Donaldson, Anti self-dual Yang--Mills connections over complex algebraic surfaces and stable vector bundles, Proc.\ London Math.\ Soc.\ (3) \textbf{50} (1985), 1--26.

\bibitem[D2]{D2} S. K. Donaldson, Infinite determinants, stable bundles and curvature, Duke Math.\ J. \textbf{54} (1987), 231--247.


\bibitem[D3]{D3} S. K. Donaldson, Some numerical results in complex differential geometry, Pure Appl.\ Math.\ Q. \textbf{5} (2009), no.~2, 571--618.


\bibitem[Fu1]{Fu1} A. Futaki, An obstruction to the existence of Einstein K\"ahler metrics, Invent.\ Math. \textbf{73} (1983), no.~3, 437--443. 

\bibitem[Fu2]{Fu2} A. Futaki, Asymptotic Chow semi-stability and integral invariants, Internat.\ J. Math. \textbf{15} (2004), no.~9, 967--979.

\bibitem[GGOW]{GGOW} A. Garg, L. Gurvits, R. Oliveira and A. Wigderson, Algorithmic and optimization aspects of Brascamp-Lieb inequalities, via operator scaling, Geom.\ Funct.\ Anal. \textbf{28} (2018), 100--145.

\bibitem[Gu]{Gu} L. Gurvits, Classical complexity and quantum entanglement, J. Comput.\ System Sci. \textbf{69} (2004), 448--484.

\bibitem[G] {G}D. Gieseker, On the moduli of vector bundles on an algebraic surface, Ann. of Math. (2) {\bf 106} (1977), no.~1, 45--60.

\bibitem[DKLR] {DKLR} M. R. Douglas; R. L. Karp; S. Lukic and R. Reinbacher,  Numerical solution to the Hermitian Yang-Mills equation on the Fermat quintic. J. High Energy Phys.  2007,  no. 12, 083, 24 pp.

\bibitem[IP]{IP} L. Ioos and L. Polterovich, Quantization of symplectic fibrations and canonical metrics, Internat.\ J. Math. \textbf{34} (2023), no.~8, Paper No. 2350043 (arXiv:2112.00419).

\bibitem[KMS]{KMS} J. Keller, J. Meyer and R. Seyyedali, Quantization of the Laplacian operator on vector bundles, I, Math.\ Ann. \textbf{366} (2016), no.~3-4, 865--907.

\bibitem[KN] {KN}G. Kempf\ and\ L. Ness, The length of vectors in representation spaces, in {\it Algebraic geometry (Proc. Summer Meeting, Univ. Copenhagen, Copenhagen, 1978)}, 233--243, Lecture Notes in Math., 732, Springer, Berlin.



\bibitem[LN]{LN} B. Lemmens and R. Nussbaum, \emph{Nonlinear Perron-Frobenius theory}, Cambridge Tracts in Math.\ 189, Cambridge Univ.\ Press, Cambridge, 2012.

\bibitem[Ma]{Ma} T. Mabuchi, An obstruction to asymptotic semistability and approximate critical metrics, Osaka J. Math. \textbf{41} (2004), no.~2, 463--472.

\bibitem[N]{N} R. D. Nussbaum, Hilbert's projective metric and iterated nonlinear maps, Mem.\ Amer.\ Math.\ Soc. \textbf{75} (1988), no.~391.

\bibitem[PS] {PS} D. H. Phong\ and\ J. Sturm, Stability, energy functionals, and K\"ahler-Einstein metrics, Comm. Anal. Geom. {\bf 11} (2003), no.~3, 565--597.

\bibitem[Sa]  {S} Y. Sano, Numerical algorithm for finding balanced metrics, Osaka J. Math. {\bf 43} (2006), no.~3, 679--688.

\bibitem[Se]  {Se} R. Seyyedali, Numerical algorithm for finding balanced metrics on vector bundles, Asian J. Math. {\bf 13} (2009), no.~3, 311--322.



\bibitem[Se2]{Se2} R. Seyyedali, Balanced metrics and Chow stability of projective bundles over K\"ahler manifolds II, J. Geom.\ Anal. \textbf{23} (2013), no.~4, 1944--1975.

\bibitem[Se3]{Se3}  R. Seyyedali, Numerical algorithm for finding balanced metrics on polarized manifolds, preprint. .

\bibitem[SH]{SH} S. Sra and R. Hosseini, Conic geometric optimization on the manifold of positive definite matrices, SIAM J. Optim. \textbf{25} (2015), 713--739.

\bibitem[T] {T}A. C. Thompson, On certain contraction mappings in a partially ordered vector space. Proc. Amer. Math. Soc. {\bf 14} (1963), 438--443.

\bibitem[Ty]{Ty} D. E. Tyler, A distribution-free $M$-estimator of multivariate scatter, Ann.\ Statist. \textbf{15} (1987), 234--251.

\bibitem[UY]{UY} K. Uhlenbeck and S.-T. Yau, On the existence of Hermitian--Yang--Mills connections in stable vector bundles, Comm.\ Pure Appl.\ Math. \textbf{39} (1986), S257--S293.


\bibitem[W] {W}X. Wang, Balance point and stability of vector bundles over a projective manifold, Math. Res. Lett. {\bf 9} (2002), no.~2-3, 393--411.


\bibitem[W2] {W2}X. Wang, Canonical metrics on stable vector bundles, Comm. Anal. Geom. {\bf 13} (2005), 253--285.





\end{thebibliography}
\end{document}